\documentclass{article}

\begin{document}


\title{Remarks on $G_{2}$-manifolds with boundary}
\author{Simon Donaldson}

\maketitle
\newcommand{\bC}{{\bf C}}
\newcommand{\bR}{{\bf R}}
\newcommand{\bZ}{{\bf Z}}
\newcommand{\hrho}{\hat{\rho}}
\newcommand{\trho}{\tilde{\rho}}
\newcommand{\vol}{{\rm vol}}
\newcommand{\dbd}{\overline{\partial}\partial}
\newcommand{\sthird}{{\textstyle \frac{1}{3}}}
\newcommand{\ssixth}{{\textstyle \frac{1}{6}}}
\newcommand{\stthird}{{\textstyle \frac{2}{3}}}
\newcommand{\sfourthree}{{\textstyle \frac{4}{3}}}
\newcommand{\squart}{{\textstyle \frac{1}{4}}}
\newcommand{\sthreetwo}{{\textstyle \frac{3}{2}}}
\newcommand{\shalf}{{\textstyle \frac{1}{2}}}
\newtheorem{prop}{Proposition}
\newtheorem{lem}{Lemma}
\newtheorem{question}{Question}
\newtheorem{cor}{Corollary}
\newtheorem{defn}{Definition}
\section{Introduction}
In this  paper we discuss special geometric structures on manifolds of dimensions 6 and 7, and the connections between these arising in the case of  a $7$-manifold with boundary. 
Our approach is largely based on a seminal paper of Hitchin, published in the {\it Journal of Differential Geometry} \cite{kn:H1} which emphasises differential forms and volume functionals (see also [7]). In dimension 6 the primary structure of interest is a Calabi-Yau structure, i.e. a Riemannian 6-manifold with holonomy contained in $SU(3)$ or, equivalently, a complex 3-manifold with a K\"ahler metric and a  holomorphic $3$-form of constant non-zero norm. In dimension 7, the primary structure of interest is a torsion free  $G_{2}$-structure {\it i.e.} a Riemannian 7-manifold with holonomy contained in $G_{2}$. But in each case there are useful variants of these where we relax the holonomy condition and consider \lq\lq closed'' $G_{2}$ and $SL(3,\bC)$ structures. In Section 2 we review the basic theory and then discuss the elementary differential geometry of a hypersurface in a $G_{2}$ manifold. The main observation is a connection between the mean curvature of the hypersurface and the intrinsic geometry of the submanifold. In Section 3 we begin with  simple remarks about the question of deforming a closed $G_{2}$-structure on a  manifold with boundary to a torsion-free structure. Then we apply some classical Riemannian geometry to the case when the boundary is constrained to have positive mean curvature. In Section 4 we consider $G_{2}$-cobordisms and discuss connections with deformations of $SL(3,\bC)$-structures  \lq\lq tamed''by a symplectic form. We explain the possible relevance of  $G_{2}$-cobordisms to questions of Torelli type for  Calabi-Yau 3-folds. In Section 5 we consider related  questions for maximal submanifolds,  which arise as dimensional reductions and adiabatic limits of the  special holonomy theory. 

The author thanks the Simons Foundation for support of this work through the  Simons Collaboration Grant \lq\lq
Special holonomy in Geometry, Analysis and Physics''. Thanks are also due to Richard Thomas and Johannes Nordstr\"om for helpful discussions. 

\section{Basics}
\subsection{Algebraic structures}

We review the special features of 3-forms in 6 and 7 dimensions. First let $V$ be a $6$-dimensional real vector space.
\begin{defn} A $3$-form $\rho\in \Lambda^{3}V^{*}$ is called definite if for each non-zero $v\in V$ the contraction $i_{v}(\rho)\in \Lambda^{2}V^{*}$ has rank $4$.
\end{defn}
The first basic fact is that if $V$ has an orientation a definite $3$-form $\rho$ defines a complex structure $I_{\rho}$ on $V$ (and reversing the orientation changes $I_{\rho}$ to $-I_{\rho}$). To see this, for $v\in V$ let $N_{v}$ be the null space of $i_{v}(\rho)$ \ {\it i.e.}
$$    N_{v}= \{ v'\in V: i_{v'} i_{v} \rho=0\}. $$
Clearly $v$ is in $N_{v}$ and $v'$ is in $N_{v}$ if and only if $v$ is in $N_{v'}$. The condition that $\rho$ is definite is that $N_{v}$ has dimension $2$, for all non-zero $v$. For each non-zero $v'$ in $N_{v}$ the form $i_{v'}(\rho)$ induces a non-degenerate symplectic form $\Omega_{v'}$ on $V/N_{v}$. We fix the orientation on $V/N_{v}$ so that $\Omega_{v'}^{2}>0$; then the map $v'\mapsto \Omega_{v'}^{2}$ defines a conformal structure on the two-dimensional vector space $N_{v}$. If we are given an orientation of $V$ we get an induced orientation on $N_{v}$, so we get a complex structure on $N_{v}$ in the usual way. Then we define $I_{\rho}(v)\in N_{v}$ using this complex structure on $N_{v}$ and it is clear that $I_{\rho}^{2}=-1$ since $N_{I_{\rho} v}= N_{v}$.  

The second basic fact is that the definite $3$-forms form a single open orbit under the action of $GL(V)$. Any such form is equivalent to the standard model 
\begin{equation} \rho_{0}= {\rm Re} (dz_{1} dz_{2}dz_{3}) \end{equation}
on $\bC^{3}$, and the complex structure defined by $\rho_{0}$ is the standard one (for the standard orientation on $\bC^{3}$). In other words, giving an oriented real $6$-dimensional vector space $V$ with a definite form $\rho$ is equivalent to giving a $3$-dimensional complex space  with a non-zero  complex form of type $(3,0)$ and the stabiliser in $GL^{+}(V)$ of $\rho$ is isomorphic to $SL(3,\bC)$. To see this choose any non-zero $v$ in $V$,  let $v'=I_{\rho}(V)$ and set $\Omega=\Omega_{v}, \Omega'= \Omega_{v'}$.
By definition these form an orthonormal pair of $2$ forms on the $4$-dimensional vector space $V/N_{v}$ with respect to the wedge-product and it is well-known that such a pair defines a complex structure on $V/N_{v}$ such that $\Omega- i \Omega'$ is   a non-zero element of the complex exterior square.  Thus there are complex co-ordinates $z_{2}, z_{3}$ on $V/N_{v}$ such that
$$   \Omega- i\Omega' = dz_{2} dz_{3}. $$
Let $z_{1}=x_{1}+ i y_{1} $ be a standard complex co-ordinate on $N_{v}$ corresponding to the basis element $v$. 
Choose a complementary sub-space $Q$ to $N_{v}$ in $V$ such that $\rho\vert_{Q}=0$---it is easy to check from the definitions that these exist. Then $z_{1}, z_{2}, z_{3}$ become  co-ordinates on $V$ and it follows from the definitions that
$$\rho= dx_{1} \wedge {\rm Re}(dz_{2} dz_{3}) - dy_{1} \wedge {\rm Im}(dz_{2}dz_{3})= {\rm Re} (dz_{1} dz_{2} dz_{3}). $$

Given the complex structure $I_{\rho}$ defined by a definite form $\rho$ and orientation we see that there is another definite form $\tilde{\rho}$ characterised by the fact that $\rho+i \tilde{\rho}$ is  of type $(3,0)$.
The volume element defined by $\rho$ is
\begin{equation} {\rm vol_{\rho}} = \frac{1}{4}\rho\wedge \trho,  \end{equation}
 and the first variation of the volume form, with respect to a variation $\delta \rho$ is
\begin{equation}  \delta {\rm vol_{\rho}}= \shalf (\delta \rho) \wedge \trho. \end{equation}
We also have the variation in $\trho$. This is given by
\begin{equation}  \delta \trho = - i (\delta\rho)_{3,0} - i (\delta\rho)_{2,1} + i (\delta\rho)_{1,2} + i (\delta\rho)_{0,3} , \end{equation}
in terms of the decomposition of $\delta\rho$ into bi-type determined by the complex structure $I_{\rho}$. 

\

Now let $U$ be a $7$-dimensional real vector space.
\begin{defn} A $3$-form $\phi\in \Lambda^{3}U^{*}$ is called definite if for each non-zero $u\in U$ the contraction $i_{u}(\phi)\in \Lambda^{2}U^{*}$ has rank $6$. 
\end{defn}
Recall that a real 2-form $\omega$ on a complex vector space is called a {\it taming form} if $\omega(\xi, I\xi)>0$ for all non-zero vectors $\xi$. This is equivalent to saying that the (1,1) component of $\omega$ is a positive $(1,1)$ form in the standard sense. For $U$ as above fix a non-zero vector $\nu\in U$ and complementary subspace $V\subset U$, so we have a fixed isomorphism $U= \bR \nu \oplus V$.  We can write any $\phi\in \Lambda^{3} U^{*}$ as
\begin{equation}  \phi= \omega\wedge dt + \rho \end{equation}
where $\omega\in \Lambda^{2}V^{*}, \rho\in \Lambda^{3}V^{*}$ and $dt\in U^{*}$ is dual to $\nu$. It is easy to check from the definitions that $\phi$ is definite if and only if $\rho$ is a definite $3$-form on $V$ and $\omega$ is a taming form for the complex structure induced on $V$ by $\rho$ and one of the orientations on $V$.  Now suppose that $U$ has a fixed orientation. We say that a definite form $\phi$ is {\it positive} if, in the description above, $\omega$ is a taming form for the complex structure defined by $\rho$ and the induced complex structure on $V$. By continuity this is the independent of the choice of vector $\nu$ and complementary subspace $V$. It is also equivalent to saying that 
\begin{equation}  (i_{u}\phi)^{2}\wedge \phi >0\end{equation}
 for all non-zero $u\in U$. Then the expression on the left hand side of (6) defines a conformal structure on $U$ and a Euclidean structure $g_{\phi}$ in this conformal class can be fixed by requiring that $\vert \phi\vert^{2}_{g_{\phi}}=7$.
The condition that $\nu$ is orthogonal to $V$ in terms of the representation (5) is that $\omega\wedge \rho=0$, or equivalently that $\omega$ has type $(1,1)$ with respect to the complex structure on $V$. The condition that $\nu$ has length $1$ is that $  \ssixth \omega^{3}= {\rm vol}_{\rho}$.   It follows that any positive form is equivalent to the model on $\bR\oplus \bC^{3}$
 $$  \phi_{0} = \omega_{0} dt + \rho_{0}, $$
where $\omega_{0}= dx_{1} dy_{1} + dx_{2} dy_{2} + dx_{3} dy_{3}$ is the standard symplectic form on $\bC^{3}$ and $\rho_{0}={\rm Re} (dz_{1}, dz_{2} dz_{3})$, as above. Thus the positive forms make up a single orbit for the action of $GL^{+}(U)$ on $\Lambda^{3}U^{*}$. 

Given a positive $3$-form $\phi$ as above we get a $4$-form $*_{\phi} \phi \in \Lambda^{4}U^{*}$, where $*_{\phi}$ is the $*$-operator defined by the Euclidean structure $g_{\phi}$.  

\subsection{Hypersurface geometry}

Now we turn to differential geometry and we recall two basic facts.
\begin{enumerate}
\item If $N$ is an oriented $6$-manifold and $\rho$ is a $3$-form on $N$ which is definite at every point then the induced almost-complex structure structure $I_{\rho}$ is integrable if and only if $d\rho=0$ and $d\tilde{\rho}=0$.
\item  If $M$ is an oriented $7$-manifold and $\phi$ is a $3$-form on $M$ which is positive at each point then $\phi$ is covariant constant with respect to the Levi-Civita connection of the Riemannian metric $g_{\phi}$ induced by $\phi$ if and only if $d\phi=0$ and 
$d*_{\phi}\phi=0$. In this case we say that $\phi$ is a torsion-free $G_{2}$-structure on $M$. 
\end{enumerate}

In this paper we want to consider relaxing these conditions so we say that:
\begin{itemize}
\item a {\it closed $SL(3,\bC)$-structure} on an oriented $6$-manifold is given by a definite $3$-form $\rho$ with $d\rho=0$;
\item a {\it closed $G_{2}$-structure} on an oriented $7$-manifold is given by a positive $3$-form $\phi$ with $d\phi=0$.
\end{itemize}

\begin{lem}
If $\rho$ is a closed $SL(3,\bC)$-structure on an oriented $6$-manifold $N$ then the $4$-form $d\trho$ has type $(2,2)$ with respect to the almost-complex structure $I_{\rho}$.\end{lem}

This is clear from Hitchin's variational point of view. Any definite form $\rho$ defines a volume form 
$\vol_{\rho}$ as above. By (3) the variation of the volume with respect to a compactly supported variation $\delta \rho$ of $\rho$ is
$$ \int_{N} \delta \rho \wedge \trho. $$
Let $v$ be a compactly-supported vector field on $N$ and $\delta\rho$ be the variation given by the Lie derivative ${\cal L}_{v}\rho$. Since $\rho$ is closed this is $ d (i_{v}\rho)$. Diffeomorphism invariance of the volume implies that
$$  \int_{N} d (i_{v}\rho) \trho =\int_{N} i_{v}\rho \wedge d\trho=0. $$
Since this is true for all $v$ we must have $i_{v}\rho \wedge d\trho=0$, pointwise for all tangent vectors $v$,  which is just the condition that $d\trho$ has type $(2,2)$.

One of the main observations in this paper is that, as a form of type $(2,2)$, there is a notion of {\it positivity} of  the tensor $d\trho$. In fact, using the volume form,  the 4-forms of type $(2,2)$ can be identified with the Hermitian forms on $T^{*}N$, for which we have the standard notion of positivity.  Any $(2,2)$-form $\sigma$  can be written in suitable co-ordinates at a point as
\begin{equation} -\squart(  \lambda_{1} dz_{2}d\overline{z}_{2}dz_{3}d\overline{z}_{3}+ \lambda_{2} dz_{3}d\overline{z}_{3}dz_{1}d\overline{z}_{1}+\lambda_{3} dz_{1}d\overline{z}_{1}dz_{2}d\overline{z}_{2}). \end{equation}  The form is positive if  each $\lambda_{i}>0$. We also consider the weaker notion of {\it semipositivity}, by which we mean that all $\lambda_{i}$ are non-negative and at least one is strictly positive. Equivalently, a  $(2,2)$ form $\sigma$ is semi-positive if $\omega \wedge \sigma>0$ for all positive $(1,1)$-forms $\omega$.   We will say that a closed $SL(3,\bC)$ structure $\rho$ is {\it mean-convex} if $d\trho$ is semi-positive at each point, and {\it strictly mean-convex} if it is positive at each point. Likewise for {\it mean-concave}. Changing the sign of the almost-complex structure interchanges the two conditions.  

  There is a scalar invariant ${\rm det}(\sigma)$ of a $(2,2)$-form $\sigma$
on a manifold $N$ with $SL(3, \bC)$ structure. To define this we use the volume
form to identify $\sigma$ with an element $\underline{\sigma}$ of $\Lambda^{2}TN$,
then take $\underline{\sigma}^{3}\in \Lambda^{6}TN$ and multiply by the $\ssixth$ of the volume
form,  to get a scalar (the factor being a convenient normalisation). In terms of the explicit representation (7) in standard
co-ordinates at a point we have $\det{\sigma}=  \lambda_{1}\lambda_{2}\lambda_{3}$. Thus for a strictly mean-convex structure $\rho$ the function $\det(d\trho)$ is strictly positive.

\

Now consider an oriented $7$-manifold $M$ with torsion-free $G_{2}$-structure $\phi$ and a $6$-dimensional submanifold $N\subset M$. By the discussion in 2.1 above the restriction of $\phi$ to $N$ is a closed definite $3$-form $\rho$.  One first observation is that the induced Riemannian measure  on $N$ coincides with $\vert {\rm vol}_{\rho}\vert$, where the choice of sign of the volume form ${\rm vol}_{\rho}$ depends on an orientation of $N$, or equivalently a co-orientation of $N\subset M$.

Fix a choice of unit normal vector field $\nu$. From standard Riemannian theory we have at each point $p\in N$ the second fundamental form $B\in s^{2} T_{p}^{*} N$.  On the other hand we have an induced  $SL(3,\bC)$
structure $\rho$ on $N$  and a $2$-form
$\omega$ on $N$ given by the contraction $i_{\nu}(\phi)$. The $2$-form $\omega$
is a positive $(1,1)$-form with respect to the almost-complex structure defined
by $\rho$. Using this almost-complex structure  we write $B$ as a sum
$B= B_{1,1} + B_{\bC}$ where $B_{1,1}$ is the real part of a Hermitian form and $B_{\bC}$ is the real part of a complex quadratic form. The component $B_{1,1}$ is somewhat analogous to the Levi form of a real hypersurface in a complex K\"ahler manifold. Using the standard identification we let $\beta_{1,1}$ be the $(1,1)$ form on $N$  corresponding to the Hermitian form $B_{1,1}$. The mean curvature $\mu$ of $N$ in $M$ is the trace of $B$ which can be written as \begin{equation} \mu=
\beta_{1,1}\wedge\omega^{2}\   ({\rm vol})^{-1}. \end{equation}
\begin{prop}
 In this situation we have $\beta_{1,1}\wedge \omega= \shalf d\trho$, so 
 $$  \mu = \shalf (d\trho\wedge\omega )\  ({\rm vol})^{-1}. $$
In particular, if the induced $SL(3,\bC)$-structure $\rho$ is mean convex then the mean curvature $\mu$ is positive with respect to the normal direction $\nu$.
\end{prop}

\

{\bf Remarks}
\begin{itemize}\item The point of the Proposition is that it relates $d\trho$, which is an intrinsic invariant of the structure on $N$,  to the mean curvature which is an extrinsic invariant of the submanifold  $N\subset M$.
\item If $\rho$ is strictly mean convex there is a stronger statement
\begin{equation}    \mu \geq \sthreetwo \det(d\trho)^{1/3}. \end{equation}
In terms of a  representation (7) in standard co-ordinates at a point this is the arithmetic-geometric mean inequality for the $\lambda_{i}$. 
\end{itemize}

The proof of Proposition 1 is a straightforward calculation which can be done in various ways. For one approach we first observe that, since a torsion-free structure agrees with the flat model to order two at a point, it suffices to consider the case when $M$ is the flat model $\bR\times \bC^{3}$ as in (2.1). We  take $N$ to be the graph of a function $f:\bC^{3}\rightarrow \bR$ with $f$ and $df$ vanishing at the origin. The induced $3$-form (pulled back to $\bC^{3}$) is given by $$\rho= \rho_{0} + df \wedge \omega_{0}. $$
Then the formula (4) for the variation of $\trho$ shows that
$$  \trho = \trho_{0} -I (df \wedge \omega_{0}) + O(\vert z\vert^{2}), $$
where $I$ acts as $+i$ on $\Lambda^{2,1}$ and $-i$ on $\Lambda^{1,2}$. So, at the origin 
$$  d\trho= -d(I df\wedge\omega_{0})= (-2i\dbd f) \wedge\omega_{0}. $$ 
On the other hand, in the familiar way the second fundamental form $B$ at the origin is given by the Hessian of $f$ and the $(1,1)$ component $\beta_{1,1}$ is $-i\dbd f$, hence the desired formula. 

We give an another approach to the calculation for the mean curvature. Let $g$ be a function of compact support on $N$. This gives a variation vector field $g\nu$ and the variation of the induced Riemannian volume of $N$ is
$$  \delta{\rm Vol} = \int_{N} g \ \mu \  {\rm vol}_{N} $$
As we observed at the beginning of this subsection, for any hypersurface in $M$ the Riemannian volume coincides with the volume computed from the induced $SL(3,\bC)$ structure.  The variation in the induced $3$-form $\rho$ is the Lie derivative $-{\cal L}_{g\nu}(\phi)$ and since $\phi$ is closed this is
$$ \delta \rho = -d ( g i_{\nu}\phi)= -d (g \omega), $$
so by (3) the variation in the volume is
$$  \delta{\rm Vol} = -\shalf\int_{N} d(g\omega) \wedge\trho = \shalf \int_{N} (g\omega)\wedge d\trho,  $$
and since this true for all $g$ we must have $\mu {\rm vol}_{N}= \shalf \omega\wedge d\trho$. This derivation has the advantage that it shows the formula for the mean curvature applies for hypersurfaces in 7-manifolds with {\it closed} $G_{2}$-structures (but the formula for $\beta_{1,1}$ in Proposition 1 then acquires an extra, trace-free, term).

For completeness we also give a formula for the component $B_{\bC}$ of the second fundamental form. Contraction of vectors with $\rho- i\trho$ defines an isomorphism $\Lambda^{0,2}= \overline{TN}$, hence $\Lambda^{1,2}= TN^{*}\otimes \overline{TN}$. Using the Hermitian metric we identify $\overline{TN}$ with $T^{*}N$ so $\Lambda^{1,2}= T^{*}N \otimes T^{*}N$. In particular we have an embedding of the symmetric tensors  $s^{2}(T^{*}N)\subset \Lambda^{1,2}$. (These are the primitive (1,2)-forms.) The component $B_{\bC}$ of the second fundamental form is the real part of  an element of $s^{2}(T^{*}N)$, so using these identifications it corresponds to a  form $\beta_{1,2}\in \Lambda^{1,2}$. Then we have
\begin{equation} d\omega= \shalf \mu \rho -\shalf \beta_{1,2}-\shalf \overline{\beta_{1,2}}. \end{equation}
We leave the verification to the reader.

\subsection{Examples}

\begin{enumerate}
\item Take $M$ to be $\bR^{7}$ with the standard flat $G_{2}$-structure and $N$ to be the $6$-sphere with outward-pointing normal.    The induced $SL(3,\bC)$ structure is strictly mean convex and the corresponding almost-complex  structure on $S^{6}$ is the standard one. More generally,  recall that a {\it nearly K\"ahler} structure on an oriented $6$-manifold $N$ is given by a closed definite $3$-form $\rho$ and  a $2$-form $\omega$ which is a positive $(1,1)$-form with respect to the almost complex structure and satisfying
$$ \ssixth\omega^{3}=\vol_{\rho}\ \ ,\ \   d\tilde{\rho}= 2\omega^{2}      \ \ , \ \  d\omega = 3\rho . $$
So the $SL(3,\bC)$-structure is mean-convex. A nearly K\"ahler-structure defines a conical torsion-free $G_{2}$-structure on $(0,\infty)\times N$ with $$\phi= r^{3}\rho + r^{2} dr \wedge \omega. $$
In particular, we recover the flat structure on $\bR^{7}$ from the $6$-sphere in this way. 

\item We say that an $SL(3,\bC)$ structure $\rho$ on $N$ is {\it tamed} by a symplectic form $\Omega$ if $\Omega$ has positive $(1,1)$ component (as in  2.1 above). Then on a compact manifold $N$, 
$$  \int_{N} d\trho \wedge \Omega = \int_{N} \trho\wedge d\Omega = 0. $$
It follows that a mean-convex $SL(3,\bC)$ structure on $N$ does not admit a  taming form. 

\item 
We consider a dimensional-reduction related to a construction of Baraglia \cite{kn:B}. Let $T^{4}=\bR^{4}/\bZ^{4}$ be the $4$-torus and identify $H^{2}(T^{4})$ with $\bR^{3,3}$, the indefinite quadratic form defined by cup-product. We also regard $\bR^{3,3}$ as the space of constant $2$-forms on $T^{4}$. Let $\Sigma\subset \bR^{3,3}$ be a \lq \lq space-like'' surface (i.e. a $2$-dimensional submanifold on which the quadratic form restricts to a Riemannian metric). Then we have a canonical $3$-form $\rho_{\Sigma}$ on $\Sigma\times T^{4}$. If we write $f:\Sigma \rightarrow \bR^{p,q}$ for the inclusion map  then, in terms of local co-ordinates $s_{1}, s_{2}$ on $\Sigma$,
\begin{equation}  \rho_{\Sigma} = - (\frac{\partial f}{\partial s_{1}} ds_{1} + \frac{\partial f}{\partial s_{2}} ds_{2}), \end{equation}
where the partial derivatives are interpreted as constant co-efficent $2$-forms on $T^{4}$. In other words $\rho_{\Sigma}=-df$, where the right hand side is interpreted as a 3-form on $\Sigma\times T^{4}$.  The submanifold $\Sigma\subset \bR^{3,3}$ has a mean curvature vector $\mu_{\Sigma}$ which is normal to $T\Sigma$. Suppose that, at a point, $\mu_{\Sigma}$ is a spacelike vector in $\bR^{3,3}$. Then $T\Sigma + \bR \mu_{\sigma} $ is a maximal positive subspace in $\bR^{3,3}$ and as such has a canonical orientation. Thus if we are given an orientation of $\Sigma$ it makes sense to say that $\mu_{\Sigma}$ is \lq\lq outward pointing''. 
\begin{lem} In  this situation 
$$  d\trho_{\Sigma} = {\rm vol}_{\Sigma}\otimes \mu_{\Sigma}, $$
where ${\rm vol}_{\Sigma}$ is the induced area $2$-form on $\Sigma$ and  $\mu_{\Sigma}$ is viewed as a constant coefficient $2$-form on $T^{4}$.
The $3$-form $\rho_{\Sigma}$ is  mean-convex   if and only if $\mu_{\Sigma}$ is space-like and outward pointing. 
\end{lem}
We leave the proof as an exercise for the  reader. Note that in this situation
$d\rho_{\Sigma}$ has rank at most $1$, so $\rho_{\Sigma}$ is never strictly mean-convex.

\

Let $\Xi \subset \bR^{3,3}$ be a 3-dimensional space-like submanifold, with inclusion map $F:\Xi\rightarrow \bR^{3,3}$. For the same reason as above it inherits an orientation, so we have an induced volume form $\chi$. We define a $3$-form on $\Xi\times T^{4}$:
$$  \phi_{\Xi}= - dF + \chi, $$
using notation as above. Then Baraglia shows that this defines a torsion-free $G_{2}$-structure if and only if $\Xi$ is a {\it maximal submanifold} (that is, a stationary point for the induced volume function, with respect to compactly supported variations in $\Xi$).  For a surface  $\Sigma\subset \Xi$, Proposition 1 amounts to the elementary statement that if $\mu_{\Sigma}$ is spacelike and outward pointing then the mean curvature of $\Sigma\subset \Xi$ is positive. 

\item We consider a different dimension reduction, as in \cite{kn:D1}. This time we take a flat $3$-torus $\bR^{3}/\bZ^{3}$ and an oriented $3$-manifold $Y$. Let $\sigma_{1},\sigma_{2}, \sigma_{3}$ be closed $2$-forms on $Y$ which are linearly independent at each point and let $\theta_{i}$ be co-ordinates on $\bR^{3}$. Then we have a definite $3$-form
$$ \rho=     d\theta_{1} d\theta_{2}d\theta_{3} - \sum \sigma_{i} d\theta_{i}$$
on the $6$-manifold $Y\times T^{3}$. By elementary linear algebra there is a unique basis of $1$-forms $\epsilon_{i}$ such that 
$$   \sigma_{i}= \epsilon_{j}\wedge \epsilon_{k}, $$
for cyclic permutations $(ijk)$ of $(123)$. One finds that
\begin{equation}  \trho = - \epsilon_{1}\epsilon_{2}\epsilon_{3}+ \sum_{{\rm cyclic}} \epsilon_{i}d\theta_{j}d\theta_{k}. \end{equation}
 The condition that $\sigma_{i}$ are closed means that
 $$  d\epsilon_{i} =\sum S_{ij} \sigma_{j}, $$
 where $(S_{ij})$ is a symmetric matrix and we have
 $$  d\trho= \sum S_{ai} \sigma_{a}\wedge d\theta_{j} d\theta_{k} $$
 (where $a$ ranges over $1,2,3$ and $(ijk)$ over cyclic permutations).
 The condition that $\phi$ is mean-convex is that $S$ is a nonnegative matrix (and not $0$). 
\end{enumerate}
\section{Manifolds with boundary}
\subsection{Gluing closed forms}
The main focus of this paper is a compact oriented $7$-manifold $M$ with boundary $N$ and a given closed $SL(3,\bC)$-structure $\rho$ on $N$. In addition we consider an \lq\lq enhancement'' of $\rho$ which is an equivalence class of closed $3$-forms over $M$ equal to $\rho$ on $N$ under the equivalence relation $\psi\sim \psi+d\alpha$ where $\alpha$ is a $2$-form vanishing on $N$. Thus the existence of an enhancement is the condition that the cohomology class $[\rho]$ extends to $H^{3}(M)$ and the difference of two enhancements is naturally an element of $H^{3}(M,N)$. We write $\hrho$ for an enhancement class. Then we have two existence questions.
\begin{itemize}
\item Is there a closed $G_{2}$-structure with boundary value $\hrho$?
\item Is there a torsion-free $G_{2}$-structure with boundary value $\hrho$?
\end{itemize}
There are also corresponding uniqueness questions. In \cite{kn:D2} we showed that the second question, modulo diffeomorphisms fixing the boundary pointwise, corresponds to an elliptic boundary value problem of index $0$. 

The next Proposition is a reflection of the fact that closed $G_{2}$-structures form a more flexible class than the torsion-free structures. Let 
 $M_{1}, M_{2}$ be oriented $7$-manifolds with boundary (not necessarily compact) and suppose that $N_{1}, N_{2}$ are  compact components of the boundary. Write $\iota_{i}:N_{i}\rightarrow M_{i}$ for the inclusion maps. 
Suppose that there is a diffeomorphism $\gamma:N_{1}\rightarrow N_{2}$ which is orientation reversing (for the orientations induced from $M_{i})$. Then in the standard way we can form a manifold $M_{1} \sharp_{\gamma} M_{2}$ by gluing the boundary components $N_{1}, N_{2}$ using $\gamma$.
\begin{prop}
Suppose that  $\phi_{1}, \phi_{2}$ are closed $G_{2}$-structures on $M_{1}, M_{2}$ and that $\iota_{1}^{*}(\phi_{1})= \gamma^{*}\iota_{2}^{*}(\phi_{2})$.
Then there is closed $G_{2}$ structure $\phi$ on $M_{1}\sharp_{\gamma} M_{2}$. Moreover $\phi$ can be chosen arbitrarily close to $\phi_{i}$ outside an arbitrarily small neighbourhood of $N_{i}$. 
\end{prop}

We sketch a proof. Assume for simplicity that $M_{i}$ are compact with just one boundary component $N_{i}$, so $M= M_{1}\sharp_{\gamma} M_{2}$ is a closed manifold. We have an $L^{\infty}$ $3$-form $\Phi_{0}$ on $M$, defined to be equal to $\phi_{i}$ on ${\rm int}\ M_{i}\subset M$. It follows from the hypothesis and Stokes' formula that $d\Phi_{0}=0$ in the weak sense. Let $K_{\epsilon}$ be the operators defining the $1$-parameter heat semigroup for the Hodge Laplacian on $3$-forms on $M$ (with some choice of Riemannian metric) and for $\epsilon>0$ set $\Phi_{\epsilon}=K_{\epsilon}\Psi_{0}$. Then $\Phi_{\epsilon}$ is a smooth, closed $3$-form and the only point to check is that $\Phi_{\epsilon}$ is positive. Near the boundary $N_{1}\subset M_{1}$ we choose a collar neighbourhood with normal co-ordinate $t\in [0,\delta)$ so that the $3$-form is
\begin{equation} \rho_{t} + \omega_{t} dt , \end{equation} 
where $\rho_{t}, \omega_{t}$ are $t$-dependent forms on $N_{1}$. The positivity condition is that $\rho_{t}$ is  definite and  that $\omega_{t}$ has positive $(1,1)$-component for the almost-complex structure defined by $\rho_{t}$. Gluing the corresponding representation of $\phi_{2}$, we can write $\Phi_{0}$ in the same form (13),  but now with $t$ in an interval $(-\delta, \delta)$ and piecewise-smooth forms $\rho_{t}, \omega_{t}$. The $3$-forms $\rho_{t}$ are continuous across $t=0$ but $\omega_{t}$ has a jump discontinuity. The key point now is that the set of $2$-forms with positive $(1,1)$ part (with respect to a fixed complex structure) forms a convex cone. In this description, the action of $K_{\epsilon}$ is given, to a very good approximation, by a positively weighted average of the forms at nearby points. It follows easily from this that
$K_{\epsilon}\Phi_{0}$ is a positive $3$-form for small $\epsilon$.

Alternatively, we can  construct suitable smoothing operators like $K_{\epsilon}$ by explicit local formulae,  thus avoiding use of any analytical theory of the heat equation. Such a construction also applies when $M_{i}$ are not compact, or have additional boundary components. Moreover we can arrange that $\phi$ is exactly equal to $\phi_{i}$ outside  arbitrarily small neighbourhoods of $N_{i}$. 

\subsection{A fundamental difficulty}

Going back to the general questions at the beginning of Section 3.1: the most optimistic, naive, hope would be that any closed $G_{2}$-structure can be deformed, through closed $G_{2}$-structures with fixed boundary data, to a torsion-free structure. There are various reasons why this cannot be true and we will discuss one such difficulty in this subsection.

   For small $\lambda>0$ let $\Omega_{\lambda}$ be a bounded domain with smooth boundary in $\bR^{7}=  \bC^{3}\times \bR$ which is diffeomeorphic to a ball and which near the origin is given by 
\begin{equation}       \{ (z,t)\in \bC^{3}\times \bR: 0< t< \vert z\vert^{2}+ \lambda \}.\end{equation}
The boundary of $\Omega_{\lambda}$ is an embedded sphere $\iota_{\lambda}: S^{7}\rightarrow \bR^{7}$ and the embeddings can clearly be taken to have a smooth limit $\iota_{0}$ which is an immersion. Let $\phi$ be the standard flat $G_{2}$-structure on $\bR^{7}$ and $\rho^{(\lambda)}=\iota_{\lambda}^{*}(\phi)$ for $\lambda\geq 0$. Then  for $\lambda>0$ the $3$-form $\rho^{(\lambda)}$ is the boundary value of a torsion-free $G_{2}$-structure on $B^{7}$ but these have no smooth limit as $\lambda$ tends to $0$. On the other hand we will show that $\rho^{(0)}$ is the boundary value of a closed $G_{2}$ structure on the ball.

Let $\pi:\bC^{3}\times \bR\rightarrow \bC^{3}\times \bR$ be the map $\pi(z,t)= (z, \vert z\vert^{2} t)$ and let $\Phi$ be the $3$-form $\Phi= \pi^{*}\phi$ on a region 
$$U_{\kappa}= \{ (z,t)\in \bC^{3}\times \bR: \vert z\vert\leq \kappa, 0\leq t\leq 1 \}. $$  Thus
$$  \Phi= \rho_{0} + t (d\vert z\vert^{2}) \omega_{0}) + \vert z\vert^{2} \omega_{0} dt. $$
We can fix a small $\kappa>0$ so that, for for this range of $t$ and $\vert z\vert$,  the $3$-form  $\rho_{0}+ t (d\vert z\vert^{2}) \omega_{0}$ is  definite on $\bC^{3}$ and  $\omega_{0}$ has positive $(1,1)$ part with respect to this form. Then $\Phi$ is a positive $3$-form on $U_{\kappa}$ except for the points where $z=0$. Let $\eta$ be a $1$-form on $\bC^{3}$ with $d\eta=\omega_{0}$ and let $\chi$ be a standard cut-off function on $\bC^{3}$, vanishing for $\vert z\vert\geq \kappa/2$ and equal to $1$ for $\vert z \vert\leq \kappa/4$.   For small $\epsilon>0$ let $$\Phi_{\epsilon}= \Phi + \epsilon d (\chi \eta dt). $$
So 
$$  \Phi_{\epsilon}= \omega_{0} (\vert z\vert^{2} +\epsilon \chi) dt + \epsilon d\chi \wedge \eta dt  + (\rho_{0}+ t (d\vert z\vert^{2}) . $$
A moments thought shows that $\Phi_{\epsilon}$ is a closed positive $3$-form on $U_{\kappa}$ for small $\epsilon$.  By construction, $\Phi_{\epsilon}$ has the same boundary values as $\Phi$ on the boundaries $t=0,1$ and agrees with $\Phi$ for $\vert z\vert>\kappa/4$. Now it is clear that we can choose a smooth map $I_{0}: B^{7}\rightarrow \bR^{7}$ extending the immersion $\iota_{0}$ and choose a region $\tilde{U}\subset B^{7}$ such that there is a diffeomorphism $h:\tilde{U}\rightarrow U_{\kappa}$ with $I_{0}= \pi\circ h$ on $\tilde{U}$ and such that $I_{0}$ is an immersion outside $\tilde{U}$.  Then the $3$-form $\phi$ which is equal to $h^{*}\Phi_{\epsilon}$ on $\tilde{U}$ and to $I^{*}(\phi_{0})$ outside $\tilde{U}$ is a closed $G_{2}$-structure on $B^{7}$ with boundary value $\rho^{(0)}$. 

This example does not completely rule out the possibility that $\rho^{(0)}$ is the boundary value of a torsion-free $G_{2}$-structure (because there could be some other solution which is not the limit of the flat solutions for $\lambda>0$), but it seems unlikely that this happens. In any case this phenomenon---of different parts of the boundary coming together---will be a serious problem in any kind of existence theory.

\subsection{Some Riemannian geometry}

In this subsection we consider a compact Riemannian manifold $X$ of dimension $(n+1)$ with smooth boundary $Y$ such that
\begin{itemize}
\item The Ricci curvature of $X$ is non-negative;
\item The mean curvature $\mu $ of the boundary (with respect to the outward pointing normal) is bounded below by a positive constant $\mu_{0}$.
\end{itemize}
We recall four results, each of a standard nature, which hold in this situation.
Let ${\cal P}$ be the set of smooth maps $\gamma:[0,1]\rightarrow X$ with $\gamma(0), \gamma(1)\in Y$ but with $\gamma(t)$ in the interior of $X$ for $0<t<1$. For $\delta>0$ we write ${\cal P}_{\delta}$ for the subset of ${\cal P}$ given by paths of length at most $\delta$.

\begin{prop} For any path $\gamma$ in ${\cal P}$ there is a small variation in ${\cal P}$ which decreases the length.
\end{prop}
By considering the first variation it suffices to consider the case when $\gamma$ is a geodesic which is normal to the boundary at the end points. Take an orthonormal frame of  $TY_{\gamma(0)}$ and parallel transport these along $\gamma$ to get variation vector fields $v_{i}$. The second variation of arc length under the variation $v_{i}$ (adapted to lie in ${\cal P}$ in the obvious way) is 
$$  -\int_{\gamma} K(v_{i}, \gamma') - B_{\gamma(0)}(v_{i}(0)) - B_{\gamma(1)}(v_{i}(1)), $$ where $B$ is the second fundamental form of the boundary and $K(\ ,\ )$ is the sectional curvature. Summing over $i$, the sum of the second variations is
$$  -\int_{\gamma} {\rm Ric}(\gamma') - \mu(\gamma(0))-\mu(\gamma(1))<0, $$ so at least one of the variations decreases length. 

\begin{prop} Let $\delta$ be the minimum length of a geodesic segment in ${\cal P}$ which is orthogonal to $Y$ at $\gamma(0)$. Then any path in ${\cal P}_{\delta}$ can be contracted to a point through paths in ${\cal P}_{\delta}$. 
\end{prop}

This follows from the previous result and an argument of Morse-theory type. 

\begin{prop}
The distance of any point of $X$ to the boundary $Y$ is at most $n \mu_{0}^{-1}$.
\end{prop}

Let $x_{0}$ be a point in the interior of $X$ and $\gamma$ minimise length among paths from $x_{0}$ to the boundary. Let $v_{i}$ be a parallel  orthonormal frame along $\gamma$ as before and consider the variation vector fields $t v_{i}$ (where we assume that $\gamma$ is parametrised by arc-length). The second variation formula shows that if the sum of the second variations is positive then the length of $\gamma$ is at most $n \mu_{0}^{-1}$.  (Equality holds when $X$ is a ball in $\bR^{n+1}$ with centre $x_{0}$).

\begin{prop}
${\rm Vol}(X)\leq \frac{n}{(n+1) \mu_{0}}  {\rm Vol}(Y)$ \end{prop}

This follows from a variant of the Bishop comparison inequality, see \cite{kn:HK}.

\

The relevance of these results for our purposes is that the hypotheses are satisfied for $(X,Y)=(M,N)$ where $M$ has a torsion-free $G_{2}$-structure and the boundary $SL(3,\bC)$ structure on $N$ is strictly mean-convex. We define  $$ m(\rho)= {\rm min}_{N} ({\rm det}\  d\trho)^{1/3} $$
so  $\mu(\rho)>0$ and by (9)  we can take $\mu_{0}= \sthreetwo m(\rho)$

\begin{itemize}
\item  Proposition 3 shows, roughly speaking, that the phenomenon discussed in 3.1 cannot occur for mean-convex boundary data. So one can be more optimistic about an existence theory for torsion-free $G_{2}$-structures with  prescribed boundary data in the case when this boundary data is mean-convex. 
\item Proposition 6 gives an upper bound 
\begin{equation} {\rm Vol}(M)\leq \frac{4}{7m(\rho)} {\rm Vol} (N) . \end{equation}
The point here is that the right hand side is entirely determined by the $SL(3,\bC)$ structure on $N$. Note that equality holds when $M$ is a ball in $\bR^{7}$.  There are reasons to expect that a torsion-free $G_{2}$-structure maximises the volume among all closed $G_{2}$ structures with given boundary data (see the discussion in \cite{kn:D2}). This raises the question whether the inequality (15) is true for {\it closed} $G_{2}$-structures on $M$, with strictly mean-convex boundary. 

\

There is a variant of this discussion for submanifolds, related to Example 3 in 2.3. Let $\Sigma\subset \bR^{p,q}$ be an oriented space-like $(p-1)$-dimensional submanifold with spacelike, outward-pointing,  mean curvature $\mu_{\Sigma}$. Suppose that $\Sigma$ is the oriented boundary of a spacelike submanifold $\Xi$. Let $\nu$ be the outward pointing unit normal to $\Sigma$ in $\Xi$. Then the mean curvature $\mu$ of $\Sigma$ in $\Xi$ is $\langle \mu_{\Sigma}, \nu\rangle$. Now we have an elementary inequality
$$  \langle \mu_{\Sigma}, \nu\rangle \geq \Vert \mu_{\Sigma}\Vert = \sqrt{\langle \mu_{\Sigma}, \mu_{\Sigma}\rangle}. $$
If $\Xi$ is a maximal submanifold the Ricci curvature of the induced metric is non-negative and we deduce from Proposition 6 that
\begin{equation} {\rm Vol}(\Xi) \leq \frac{p-1}{p} {\rm Vol}(\Sigma) \left(\min_{\Sigma} \Vert \mu_{\Sigma}\Vert\right)^{-1}. \end{equation}
(With equality for a standard ball in $\bR^{p}\subset \bR^{p,q}$.) The question that arises is whether this holds for {\it any} spacelike $\Xi$ with boundary $\Sigma$. 
\end{itemize}

\section{$G_{2}$-cobordisms}

In this section we consider a pair of compact, connected, 6-manifolds with closed $SL(3,\bC)$ structures $(N_{0}, \rho_{0}), (N_{1}, \rho{1})$ and a cobordism $M$ from $N_{0}$ to $N_{1}$ with closed or torsion-free $G_{2}$ structure $\phi$ restricting to $\rho_{i}$ on the boundary. Proposition 2 shows that the existence of a closed $G_{2}$-cobordism defines a transitive relation on $SL(3,\bC)$ structures, but the orientations in the set-up mean that this relation is  not symmetric (or at least, not in an obvious way).

Our main focus is on the case when $N_{0},N_{1}$ are diffeomorphic (so we just write $N$) and $M$ is a product, as a smooth manifold. Choosing such a product structure we can express a $3$-form $\phi$ in the usual way as $\rho_{t} + \omega_{t} dt$. The existence of a closed $G_{2}$-cobordism is equivalent to the existence of a path $\rho_{t}$ from $\rho_{0}$ to $\rho_{1}$ through closed $SL(3,\bC)$-structures on $N$ such that 
\begin{equation}    \frac{d\rho_{t}}{dt}= d \omega_{t}\end{equation} where $\omega_{t}$ has positive $(1,1)$ part with respect to $\rho_{t}$. Of course we need to assume that $\rho_{0}, \rho_{1}$ define the same cohomology class in $H^{3}(N)$.  If we ignore the positivity condition  it is known that we can find  some path $\rho_{t}$ of $SL(3,\bC)$-structures. In \cite{kn:CN}, Crowley and Nordstr\"om show that \lq\lq co-closed'' $G_{2}$-structures on $7$-manifolds obey an h-principle and the same arguments apply to closed $SL(3,\bC)$-structures \cite{kn:O}. (The main point is that any hypersurface in $\bR^{7}$ acquires a closed $SL(3,\bC)$-structure,  just as any hypersurface in $\bR^{8}$ acquires a co-closed $G_{2}$-structure.)
 Easy bundle  theory considerations show that $\rho_{0}, \rho_{1}$ are homotopic as definite $3$-forms and the h-principle shows that these forms can be taken to be closed, in a fixed cohomology class. 

\

{\bf Example} Consider the standard closed definite form $\rho$ on $S^{6}$. Then $-\rho$ is also a closed definite form and there is an obvious homotopy $$ \rho_{t}= \cos(\pi t) \rho + \sin (\pi t) \trho, $$
through definite forms, but these are not closed. The Crowley-Nordstr\"om theory shows that there is some homotopy through closed definite forms. Note that such a homotopy cannot be invariant under $G_{2}$, acting on $S^{6}$. It is interesting to ask whether there is a closed $G_{2}$-cobordism from $\rho$ to $-\rho$ (or from $-\rho$ to $\rho$).

\subsection{Taming forms and cobordisms}

There is a further connection between homotopy of definite forms and $G_{2}$-cobordism in the presence of a taming form. 
\begin{lem} Let $\Omega$ be a symplectic form on $N$. If $\rho_{0}, \rho_{1}$ can be joined by a path $\rho_{t}$ of closed definite $3$-forms in a fixed cohomology class such that $\rho_{t}$ is tamed by $\Omega$ for each $t$ then there is a closed $G_{2}$-cobordism from $\rho_{0}$ to $\rho$.
\end{lem}

The proof is easy:
the hypotheses mean that we can find closed, $\Omega$-tamed, $SL(3,\bC)$ structures $\rho_{t}$ and $2$-forms  $\tilde{\omega}_{t}$ such that $d \tilde{\omega}_{t}$ is the $t$-derivative of $\rho_{t}$. Then we set $\omega_{t}=\tilde{\omega}_{t} + A \Omega$ for large $A$.

We  say that a closed $G_{2}$-cobordism is tamed by $\Omega$ if there is a product structure $M=N\times [0,1]$ with respect to which $\rho_{t}$ is  tamed by $\Omega$ for all $t$. 

There is a more precise statement of this Lemma, involving the \lq\lq enhancement'' of the boundary data. For simplicity we consider the case when $H^{2}(N)=\bR$ and fix a pair of $2$-cycles $S_{0}, S_{1}$ representing a generator of $H_{2}(N)$. Also fix a 3-chain $W\subset M= N\times [0,1]$ with boundary $-S_{0}$ in one end and $S_{1}$ in the other. So for any closed form $\phi$ with boundary values $\rho_{0}, \rho_{1}$ we have a real number 
$$    I_{W}(\phi) =\int_{W} \phi, $$
which depends only on the relative homology class of $W$. The more refined question is to ask for which values of $I_{W}$ (if any) is there a closed $G_{2}$-cobordism
from $\rho_{0}$ to $\rho_{1}$. In the presence of a symplectic form $\Omega$ as above, fix the sign of $S_{1}$ so that $\langle \Omega, S_{1}\rangle >0$.  Then
$$  \int_{W} \Omega dt = \int_{W} d(t\Omega) = \int_{S_{1}} \Omega >0.$$ The more precise statement of Lemma 3 is that, under the hypotheses of the Lemma,  there is a $\kappa_{0}$ such that for all $\kappa\geq \kappa_{0}$ there is a closed $G_{2}$-cobordism $\phi$ from $\rho_{0}$ to $\rho_{1}$ with $I_{W}(\phi)=\kappa$. Motivated by this we can define an invariant $D_{W}(\rho_{0},\rho_{1})$ of a pair $\rho_{0}, \rho_{1}$ to be the infimum of the values of $I_{W}(\phi)$ such that there is a closed $G_{2}$-cobordism $\phi$ from $\rho_{0}$ to $\rho_{1}$ (and $+\infty$ if this set is empty).  Of course this depends on the choice of $W$, but the invariants given by different choices are related  by the addition of a constant determined by  homological considerations.

There is a potential connection between these ideas and the enumerative geometry of holomorphic curves in $N$. The appropriate theory would probably be an extension  of the \lq\lq Donaldson-Thomas'' invariants to the symplectic case and, since such a theory has not so far been set up rigourously, we only sketch the idea. Suppose, in the simplest situation, that there is a single holomorphic curve for the almost complex structure defined by $\rho_{0}$ in the homology class  $[S_{0}]$ so we take $S_{0}$ to be this holomorphic curve. Similarly suppose that there is single holomorphic curve in this homology class for the almost complex structure defined by $\rho_{1}$ and take $S_{1}$ to be that curve. Suppose further that throughout the the $1$-parameter family $\rho_{t}$ there is just a single holomorphic curve $S_{t}$, giving a smoothly varying family from $S_{0}$ to $S_{1}$. This family defines a cycle $W$ in $N\times [0,1]$ and we have 
\begin{equation} I_{W}(\phi) = \int_{0}^{1} \left(\int_{S_{t}}\omega_{t} \right)  \ dt .\end{equation} 
The derivation of this equation uses the fact that for any tangent vector $v$ to $N$ at a point of $S_{t}$ the contraction $i_{v}(\rho_{t})$ vanishes on the tangent space of $S_{t}$, which is a characterisation of holomorphic curves in this setting. The point now is that $I_{W}(\phi)>0$ since $\omega_{t}$ has positive $(1,1)$ component and its integral over a curve is positive. 

Of course the situation above cannot be expected to hold in general. However we do have compactness results for holomorphic curves in the case of tamed structures and the discussion can be extended. For example we might have a finite number of holomorphic curves, with respect to $\rho_{i}$, in the given homology class and we then take the cycles $S_{i}$ to be sum of these, with suitable signs. But we will not try to go into further details here. The general point is that we can hope that there are preferred chains $W$ which impose a constraint  $I_{W}(\phi)>0$,  at least for tamed, closed, $G_{2}$-cobordisms.

\subsection{More Riemannian geometry and questions of Torelli type}

Now consider a compact Riemannian manifold with boundary which gives a cobordism from $Y_{0}$ to $Y_{1}$. Then we have another result of a standard nature.
\begin{prop}
If the Ricci curvature of $X$ is non-negative and the mean curvature of both boundary components (with respect to the outward normals) is non-negative then $X$ is isometric to a product $Y\times [0,L]$ and in particular $Y_{1}, Y_{2}$ are isometric.
\end{prop}
To prove this we consider the harmonic function $h$ on $X$, equal to $0$ on $Y_{0}$ and to $1$ on $Y_{1}$. Then we have
 \begin{equation} \Delta (\vert \nabla h\vert^{2})= \vert \nabla \nabla h\vert^{2} + {\rm Ric}(\nabla h) \end{equation}
and integrating we obtain
\begin{equation} \int_{Y_{0}} \nabla_{\nu} \vert \nabla h\vert^{2} + \int_{Y_{1}} \nabla_{\nu} \vert \nabla h\vert^{2} + \int_{X} \vert \nabla \nabla h\vert^{2} + {\rm Ric} (\nabla h) =0, \end{equation}
where $\nabla_{\nu}$ denotes the normal derivative. The second fundamental form of the boundary is the quadratic form defined by  $B(\xi)= \langle \nabla_{\xi} \nu, \xi\rangle$ for vectors $\xi$ tangent to the boundary. So if $\xi_{i}$ is an orthornomal frame the mean curvature is
$$  \mu= \sum_{i} \langle \nabla_{\xi_{i}}\nu, \xi_{i}\rangle. $$
On the boundary, write $\nabla h= f \nu$. Then 
$$  \sum_{i} \nabla_{\xi_{i}}\langle \nabla h, \xi_{i}\rangle= \sum \langle \nabla_{\xi_{i}}f \nu, \xi_{i}f \rangle = f^{2} \mu. $$
On the other hand 
$$  \Delta h = \sum_{i} \langle \nabla_{\xi_{i}} \nabla h, \xi_{i}\rangle + \langle \nabla_{\nu}\nabla h, \nu\rangle, $$
so $ f \mu= - \langle \nabla_{\nu} \nabla h, \nu\rangle$, since $h$ is harmonic. Thus $$\nabla_{\nu} (\vert \nabla h \vert^{2}) = - f^{2} \mu= -\mu \vert \nabla h\vert^{2}$$ and (20) becomes
$$     \int_{Y_{0}} \mu \vert \nabla h \vert^{2} + \int_{Y_{1}} \mu \vert \nabla h \vert^{2} + \int_{X} \vert \nabla \nabla h \vert^{2} + {\rm Ric}(\nabla h) =0. $$
Under our hypotheses all terms are non-negative so must vanish identically. In particular $\nabla\nabla h=0$ and this leads easily the product decomposition. 
\

We apply this to the case of a torsion-free $G_{2}$-cobordism. 
\begin{cor} Let $(N_{0}, \rho_{0}), (N_{1}, \rho_{1})$ be a pair of compact 6-manifolds with {\it integrable} $SL(3,\bC)$ structures. If there is a torsion-free $G_{2}$-cobordism $M$ from $(N_{0}, \rho_{0})$ to $(N_{1}, \rho_{1})$ then they are isomorphic.
\end{cor}
This follows from  Proposition 7 since the Ricci curvature of $M$ and mean curvature of the boundary vanish (the latter by Proposition 1). 

This corollary potentially has some bearing on  questions of Torelli type for Calabi-Yau 3-folds. That is, the question whether a Calabi-Yau structure is uniquely determined by the cohomology class $[\rho]\in H^{3}(N;\bR)$. In fact the usual algebraic geometry formulation is in terms of the larger data
$[\rho+ i \trho]\in H^{3}(N,\bC)$. There are examples showing that \lq\lq global Torelli'' fails,  in the standard algebraic geometry formulation \cite{kn:Sz}. But it is possible that there could be alternative formulations with  positive answers. 

\begin{question}
Suppose $(N,\rho_{0})$ and $(N,\rho_{2})$ are integrable $SL(3,\bC)$ structures.
\begin{itemize} \item If there is a closed $G_{2}$-cobordism between the structures are they isomorphic?
\item If  $\rho_{0}, \rho_{1}$ are homotopic through tamed, closed, $SL(3,\bC)$-structures in a fixed cohomology class are they isomorphic?
\end{itemize}\end{question}

In other words, it is possible that examples where the Torelli property fails come from Calabi-Yau structures in different connected components under tamed deformations (although by the Crowley-Nordstr\"om theory discussed above they  lie in the same connected component of closed $SL(3,\bC)$ structures). To explain the relevance of Corollary 1, suppose that  $\rho_{s}$ is a path of tamed, closed, structures from $\rho_{0}$ to $\rho_{1}$. We showed in \cite{kn:D2} that for small $s$ there is a torsion-free $G_{2}$-cobordism from $\rho_{0}$ to $\rho_{s}$.  If this can be continued all the way to $s=1$ we would deduce from Corollary 1 that $\rho_{0}$ and $\rho_{1}$ are isomorphic.  
 
 Continuing in a speculative vein,  similar ideas could possibly be relevant to proving {\it existence} of Calabi-Yau structures. Suppose that $\rho_{0} $ is a real-analytic,  closed $SL(3,\bC)$ structure on $N$.  Then it is straightforward to show that there a torsion-free $G_{2}$-cobordism from $\rho_{0}$ to  some $\rho_{1}$ close to $\rho_{0}$. Fix these boundary values $\rho_{0},\rho_{1}$ and attempt to vary the enhancement data, i.e. seek torsion-free $G_{2}$-cobordisms $\phi_{L}$ with  $I_{W}(\phi)= L$ and with $L\rightarrow \infty$. The simplest picture of what could happen is that,  for a suitable family of base points in $M$,  the based Gromov-Hausdorff limit as $L\rightarrow\infty$ is a product $N\times \bR$,  with a Calabi-Yau structure on $N$.  

For another question, let $\rho_{0}$ be mean-concave and $\rho_{1}$ be mean-convex. Then Proposition 7 shows that there is no torsion-free $G_{2}$-cobordism for $\rho_{0}$ to $\rho_{1}$. (The signs  are confusing here---the condition that $\rho_{0}$ is mean concave says that the boundary has positive mean curvature with respect to the outward normal, due to the switch in orientation). This can also be seen using geodesics and the second variation formula, as in Proposition 3. The question arises whether there can be a {\it closed} $G_{2}$-cobordism from $\rho_{0}$ to $\rho_{1}$.

\section{Variants for maximal submanifolds}

We can develop the same ideas in the direction of existence and uniqueness questions for maximal submanifolds. This is related to the $G_{2}$-discussion via the dimension reduction procedure described in Example 3 of 2.3, but can also be pursued independently. Consider a space-like $(p-1)$-dimensional submanifold $\Sigma\subset \bR^{p,q}$ as in 3.2. Suppose that $\Xi_{0}, \Xi_{1}$ are two $p$-dimensional maximal spacelike submanifolds with boundary $\Sigma$. For $L>0$ we consider  $[0,L]\times \bR^{p,q}\subset \bR^{p+1, q}$ and the set
$$  T= [0,L]\times \Sigma \cup \{0\}\times  \Xi_{0}\cup \{L\} \times \Xi_{1}. $$

\begin{prop} If there is a compact $(p+1)$-dimensional maximal submanifold $Z$ in $\bR^{p+1,q}$ with boundary $T$ then $\Xi_{0}=\Xi_{1}$. 
\end{prop}

(More precisely, $Z$ should be a \lq\lq manifold with corners''.) 
To see this we follows the proof of Proposition 7. The linear projection to the first factor is a harmonic function $h$ on $Z$. The maximal condition implies that the Ricci curvature of $Z$ is nonnegative and the fact that $\Xi_{i}$ are maximal implies that the mean curvature of $\Xi_{i}$ in $Z$ vanishes. The new feature is that $Z$ has an extra boundary component $[0,L]\times \Sigma$. Let $e\in \bR^{p+1,q}$ be the co-ordinate vector corresponding to the $[0,L]$ factor. Then $\vert \nabla h\vert$ at a point of $Z$ is  the length of the orthogonal projection of $e$  to the tangent space of $Z$ (with respect to the indefinite form). At points of the boundary component $[0,L]\times \Sigma$ the first vector $e$ lies inside this tangent space so $\vert \nabla h\vert=1$. A moments thought shows that the normal derivative of $\vert \nabla h\vert^{2}$ vanishes, thus we do not get any contribution to the boundary term and the same argument applies to show that $\nabla\nabla h=0$. This means that $\vert \nabla h\vert=1$ everywhere which can only happen if $e$ is tangent to $Z$ at each point and we deduce that $\Xi_{0}=\Xi_{1}$ and $Z=[0,L]\times \Xi_{0}$.

\

Finally we consider a more complicated geometric set-up, following \cite{kn:D0} (to which we refer for more details).  Let $P$ be a $p$-dimensional manifold and $Q\subset P$ a co-oriented submanifold of co-dimension $2$. We regard $(P,L)$ as an orbifold, so we have  orbifold charts around  points of $Q$ modelled on $\bR^{p-2}\times \bC$,  with the involution $z\mapsto -z$ on the $\bC$ factor. We consider a flat affine orbifold bundle $V\rightarrow P$ with structure group the affine extension $\Gamma$ of $O(p,q)$. Thus over $P\setminus L$ we have a flat $\Gamma$-bundle in the usual sense and the orbifold structure over a point $x$ of $Q$ is given by an element $r_{x}$ of order $2$ in $\Gamma$. We suppose that the $r_{x}$ are reflections in \lq\lq timelike'' vectors. Given this data,  we have a notion of a {\it branched  section} $u$ of $V$.  By definition this is given over $P\setminus Q$ by a section of the flat bundle. Locally, over  small open sets $\Pi\subset P\setminus Q$  this is represented by a map $u_{\Pi}:\Pi\rightarrow \bR^{p,q}$ and we require that this be an embedding with image a space-like submanifold. Around a point $x$ of $Q$ the behaviour of $u$ can be described as follows. There is an orthogonal decomposition
$$  \bR^{p,q}= \bC\times \bR^{p-2}\times \bR \times \bR^{q-1}$$
in which the reflection $r_{x}$ acts as $-1$ on the $\bR$ factor and $+1$ on the other factors. The factor $\bC\times \bR^{p-2}$ is a positive subspace for the indefinite form and the factor $\bR\times \bR^{q-1}$ is a negative subspace.  We can choose local co-ordinates $(w,\tau)\in \bC\times \bR^{p-2}$ on $P$ such that $Q$ is defined by $w=0$ and the section is give by a multi-valued function
\begin{equation}    u(w,\tau)= (w,\tau, f(w^{1/2}, \tau), g(w,\tau)), \end{equation}
where $f$ is an odd function in the $w^{1/2}$ variable. In other words, the orbifold co-ordinate $z$ is $w^{1/2}$ and $f$ is a genuine function $f(z,\tau)$ with $f(-z,\tau)=-f(z,\tau)$. We require that
\begin{equation}   g=O(\vert w\vert^{2}), \nabla g = O(\vert w\vert), \nabla^{2} g = O(1),\end{equation}\begin{equation} f=O(\vert w\vert^{3/2}), \nabla f=O(\vert w\vert^{1/2}), \nabla^{2} f =0(\vert w\vert^{-1/2}). \end{equation}

Finally we can define {\it maximal branched sections} of $V$ to be branched sections which away from $Q$ are locally given by parametrised maximal submanifolds of $\bR^{p,q}$. Around points of $Q$ they correspond to branched maximal subvarities, with co-dimension 2 singularities.

 Maximal branched sections are certainly not unique.  We define an equivalence relation on branched sections of $V$ as follows. If $f:(P,Q)\rightarrow (P,Q)$ is a diffeomorphism and if there is an isomorphism $\tilde{f}: f^{*}(V)\rightarrow V$ then $u\sim \tilde{f}(f^{*}(u))$. (In particular, if $f$ is isotopic to the identity the flat structure defines a lift $\tilde{f}$.) Then if $u_{0}$ is a maximal branched section and if $u_{1}\sim u_{0}$ then so also is $u_{1}$. Locally, this just corresponds to different choices of parametrisation of the same maximal submanifold. 
  Another simple way in which uniqueness can fail occurs when there is a {\it covariant constant} section $s$  of the flat orbifold vector bundle $\underline{V}$ associated to the affine bundle $V$. In that case we can change a maximal branched section $u$ to another $u+ s$. Locally this just corresponds to translation of the maximal subvariety. In most cases of interest there will be no such covariant constant sections.

\

 Now  consider the product $P\times [0,L]$ with projection $\pi:P\times [0,L]\rightarrow P$. The pull-back $\pi^{*}(V)$ is a flat affine orbifold bundle over $(P\times[0,L], Q\times [0,L])$. We consider the bundle $\pi^{*}(V)\times \bR$ over $(P\times [0,L], Q\times [0,L])$ with the obvious structure of an affine orbifold bundle with fibre $\bR^{p+1, q}$. Let $e$ be the covariant constant section of the vector bundle $\underline{\pi^{*}(V)\times \bR}$ corresponding to the unit vector in the $\bR$ factor. If we have two branched sections $u_{0}, u_{1}$ of $V$  we can consider branched sections $U$ of $\pi^{*}(V)\oplus \bR$ with boundary conditions that $U=u_{0}$ over $P\times \{0\}$ and  $U= u_{1}+ L e $ over $P\times \{L\}$.
\begin{prop} If $u_{0}$ and $u_{1}$ are two branched maximal section of $V$ and if there is a branched maximal section $U$ of $\pi^{*}(V)\times \bR$ with these boundary values then $ u_{1}\sim u_{0}+s$ for a covariant constant section $s$ of $V_{0}$. 
\end{prop}

The section $U$ induces a  Riemannian metric $\Gamma$ on $P\times [0,L]$ with a singularity along $Q\times [0,L]$. In the local-co-ordinates given by (21) this metric is uniformly equivalent to the Euclidean metric, with Lipschitz metric tensor. As before the metric has non-negative Ricci curvature away from the singular set. We write $h$ for the function on $P\times [0,L]$ given by projection of $U$ to to the $\bR$ factor in $\pi^{V}\times \bR$. Thus $h=0,L$ on the two boundary components. The local geometry away from the singular set is just as before but we need to check that the singularity does not affect the argument. Let $N_{\epsilon}$ be a tubular neighbourhood of $Q$ of radius $\epsilon$ and consider 
$$ \int_{(P\setminus N_{\epsilon}) \times [0,L]} \Delta \vert \nabla h\vert^{2} $$
There is a new boundary term 
\begin{equation}  \int_{\partial N_{\epsilon} \times [0,L]} \nabla_{\nu} \vert \nabla h \vert^{2}. \end{equation}
The integrand is locally $\langle B(\nabla h,\nu), e\rangle$ where $B(\ ,\ )$ is the second fundamental form of the image of $U$, regarded now as a bilinear form on the tangent space. It follows from (21) that $B$ is $O(\epsilon^{-1/2})$ and $\vert \nabla h \vert $ is $O(1)$. The local bi-Lipschitz property implies that  the volume of $\partial N_{\epsilon}\times [0,L]$ is $O(\epsilon)$ so the integral in (24) is $O(\epsilon^{1/2})$ and taking $\epsilon\rightarrow 0$ we deduce that $\nabla\nabla h=0$, as before. 
 In particular the length $\vert \nabla h\vert$ is a constant $c$ and $c\geq 1$ (since it is the projection of a unit vector to a maximal positive subspace in $\bR^{p+1,q}$). The local representation (21) shows that $\nabla h$, regarded as the gradient vector field, is Lipschitz on $P\times [0,L]$ and this implies that the integral curves run from one boundary component to the other (as in the smooth case). The same argument shows that the (singular) Riemannian manifold $(P\times [0,L], \Gamma)$ is isometric to a Riemannian product, say $(P, g)\times [0,L/c]$,  where the function $h$ on $P\times [0,L]$ goes over to the function $\tilde{h}(x,t) = ct$ on $(P,g)\times [0,L/c]$.

 Suppose first that $c=1$. This implies that at each point the gradient vector $\nabla h$ in the tangent space of $P\times [0,L]$ maps under the derivative of $U$ to the fixed vector $e$. Let $F:(P,g)\times [0,L]\rightarrow P\times [0,L]$ be the diffeomorphism given by the Riemannian product structure, equal to the identity on $P\times \{0\}$.  The pull back by $F$ of $\nabla h$ is the unit vector field $\partial_{t}$ in the $[0,L]$ factor. Thus the derivative  $\partial_{t}(F^{*}(U))$ is equal to $e$. The flat structure, and the fact that $F$ is the identity on $P\times\{0\}$,  gives a canonical isomorphism $\tilde{F}: F^{*}(\pi^{*}(V)\times \bR)\rightarrow \pi^{*}V\times \bR$. So we can regard $F^{*}(U)$ as a $1$-parameter family $\left(F^{*}(U)\right)_{t}$ of sections of the bundle $V\times\bR\rightarrow P$ and our identification of the $t$-derivative shows that  
\begin{equation}   \left(F^{*}(U)\right)_{t} = u_{0} + t e . \end{equation}
Now let $F$ be given on the other boundary component by $F(x,L)= (f(x), L)$ for a diffeomorphism $f:P\rightarrow P$ and let $\tilde{f}$ be the restriction of $\tilde{F}$. Then (25) specialises to
$   \tilde{f}( f^{*}(u_{1}))= u_{0}$, which shows that $u_{1}\sim u_{0}$. 
The argument above is essentially the same as that in the proof of Proposition 8, once we know that $\vert\nabla h\vert=1$. The extra difficulty that arises now is to analyse the case when $c>1$. 
 To handle this we need a lemma from local
differential geometry.

\begin{lem}
Suppose $X$ is a connected p-dimensional Riemannian manifold (not necessarily complete) and suppose that $f_{s}:X\rightarrow \bR^{p,q}$ is a smooth family of spacelike embeddings for $s\in (-\delta,\delta)$. Fix $c>1$ and let $\bR^{p,q}\times \bR$ have the standard indefinite form, positive on the $\bR$ factor. Let
$$\Phi:X\times (-\delta,\delta)\rightarrow \bR^{p+1,q}$$
be the map $\Phi(x,s)= (f_{s}(x), cs)$. Suppose that 
\begin{enumerate}
\item $\Phi$ is an isometric embedding, with space-like image, for the Riemannian product metric on $X\times (-\delta,\delta)$;
\item the image of $\Phi$ is a maximal space-like submanifold
in $\bR^{p,q}\times \bR$;
\end{enumerate}
Then there is a vector $\nu\in \bR^{p,q}$ with $\vert \nu\vert^{2} = 1-c^{2}< 0$ such that $f_{s}(x)= f_{0}(x)+ s \nu$ and  the image $f_{0}(X)$ lies in a hyperplane normal to $\nu$.
\end{lem}

Write $\frac{\partial f_{s}}{\partial s}=\nu_{s,x}$ so $\nu$ takes values in $\bR^{p,q}$. The isometric embedding condition in item (1) is equivalent to
\begin{itemize}\item $\vert \nu_{s,x}\vert^{2}=1-c^{2}$;
\item $\nu_{s,x}$ is orthogonal to the tangent space of $ f_{s}(X)$ at $f_{s}(x)$;
\item each $f_{s}$ is an isometric embedding of $X$ in $\bR^{p,q}$.
\end{itemize}

When the codimension,  $q$,  is large these conditions admit many solutions so we have to bring in the second hypothesis, that the image of $\Phi$ is a maximal submanifold. Let $\Gamma$ be the Gauss map of the image of $\Phi$. The maximal submanifold condition implies that
$$  \vert \frac{\partial\Gamma}{\partial s}\vert^{2}= {\rm Ric}_{X\times (-\delta,\delta)}(\partial s)$$
where on the left hand side we use the standard Riemannian metric on the Grassmann manifold of maximal positive subspaces. Since, for the product manifold, this component of the Ricci curvature is zero we deduce that $\Gamma$ is constant in $s$. By simple linear algebra and the orthogonality condition this implies that $\nu_{s,x}$ is independent of $s$ so we can write $\nu_{x}$v and $  f_{s}= f_{0} + s \nu_{x}$. From this one deduces easily that $\nu_{x}$ is independent of $x$, and the orthogonality shows that $X_{0}$ lies in a hyperplane normal to $\nu$.

Given this Lemma it is easy to extend the proof that we gave for the case $c=1$ to the general case.

\

There are two situations in which this result interacts with $G_{2}$-geometry.

\begin{enumerate}
\item Take $p=2, q=19$ and $P=S^{2}$. Consider a polarised Calabi-Yau threefold $N$ which admits a holomorphic Lesfchetz fibration $N\rightarrow S^{2}$ with $K3$ fibres. The cohomology of the fibres orthogonal to the K\"ahler class defines a flat orbifold vector bundle (with $Q$ the finite set of critical values) and a class in $H^{3}(N)$ yields a lift to an affine bundle $V$. The period map of the complex structure defines a branched maximal  section. The uniqueness question is a version of the Torelli problem for $K3$-fibred Calabi-Yau 3-folds. In the argument, the maximal section $U$ over $S^{2}\times [0,L]$ corresponds to the \lq\lq adiabatic limit'' of a $G_{2}$-cobordism with a Kovalev-Lefschetz fibration, as discussed in \cite{kn:D0}.
\item  Take $p=3, q=19$. Then $P$ is a $3$-manifold and $Q\subset P $ is a link. The uniqueness question for maximal sections is the adiabatic limit of a \lq\lq Torelli problem'' for closed $G_{2}$-manifolds with Kovalev-Lefschetz fibrations.  

\end{enumerate}

\end{document}